\documentclass[11pt,a4paper]{scrartcl}

\KOMAoptions{DIV=12,bibliography=oldstyle}
\usepackage[T1]{fontenc}
\usepackage{lmodern}
\usepackage[english]{babel}
\usepackage{amsmath,amssymb,amsthm,amsfonts}
\usepackage{tikz,pgf}
\usepackage{hyperref}
\hypersetup{
    colorlinks,
    linkcolor={red!50!black},
    citecolor={blue!50!black},
    urlcolor={blue!80!black}
}

\tikzstyle{vertex}=[circle, draw, fill=black, minimum size=5.5pt, inner sep=0.8pt]
\newcommand{\vertex}{\node[vertex]}
\tikzstyle{edge}=[very thick,font=\bfseries]

\newcommand{\R}{\mathbb{R}}
\newcommand{\C}{\mathbb{C}}

\newcommand{\edg}{\mcal{E}}
\newcommand{\Lam}{\operatorname{Lam}}
\newcommand{\leftquot}[2]{{}^{#1}\!{#2}}
\newcommand{\rightquot}[2]{{#2}^{#1}}

\newcommand{\subtract}{\mathop{\backslash}}
\newcommand{\quotient}{\mathop{/}}

\newcommand{\vectwo}[2]{\binom{#1}{#2}}
\newcommand{\vecxy}{\vectwo{x}{y}}
\newcommand{\mcal}{\mathcal}

\theoremstyle{plain}
\newtheorem{lemma}{Lemma}[section]
\newtheorem{proposition}[lemma]{\textbf{Proposition}}
\newtheorem{theorem}[lemma]{\textbf{Theorem}}
\theoremstyle{definition}
\newtheorem{definition}[lemma]{\textbf{Definition}}

\title{Computing the number of realizations of a Laman graph}
\author{%
Jose Capco%
	\thanks{Research Institute for Symbolic Computation (RISC), Johannes Kepler University Linz}\, 
	\thanks{Supported by the Austrian Science Fund (FWF): P28349}\,
	\thanks{Supported by the Austrian Science Fund (FWF): W1214-N15, project DK9}
\and
Matteo Gallet%
	\thanks{Johann Radon Institute for Computational and Applied Mathematics (RICAM), Austrian Academy of Sciences}\,
	\footnotemark[3]\,
	\thanks{Supported by the Austrian Science Fund (FWF): P26607}
\and
Georg Grasegger%
	\footnotemark[4]
\and
Christoph Koutschan%
	\footnotemark[4]\,
	\footnotemark[3]\,
	\thanks{Supported by the Austrian Science Fund (FWF): P29467}
\and
Niels Lubbes%
	\footnotemark[4]\,
	\footnotemark[5]
\and
Josef Schicho%
	\footnotemark[1]\,
	\footnotemark[5]
}
\date{\vspace{-0.5cm}}

\begin{document}
\maketitle
\begin{abstract}
	Laman graphs model planar frameworks which are rigid for a general choice 
	of distances between the vertices. There are finitely many ways, up to 
	isometries, to realize a Laman graph in the plane. 
	In a recent paper we provide a recursion formula for this 
	number of realizations using ideas from algebraic and tropical geometry.
	Here, we present a concise summary of this result focusing on the main ideas and the combinatorial point of view.
\end{abstract}

\section{Introduction}
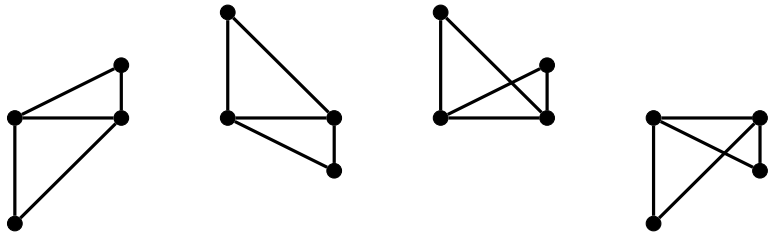
\begin{figure}[ht]
	\begin{center}
		\begin{tikzpicture}[scale=0.7]
			\begin{scope}
			  \vertex (a1) at (0,-2) {};
				\vertex (b1) at (0,0) {};
				\vertex (c1) at (2,0) {};
				\vertex (d1) at (2,1) {};
				\path[edge] (a1)edge(b1) (b1)edge(c1) (c1)edge(a1) (b1)edge(d1) (c1)edge(d1);
			\end{scope}
			\begin{scope}[xshift=1cm]
			  \vertex (a2) at (3,2) {};
				\vertex (b2) at (3,0) {};
				\vertex (c2) at (5,0) {};
				\vertex (d2) at (5,-1) {};
				\path[edge] (a2)edge(b2) (b2)edge(c2) (c2)edge(a2) (b2)edge(d2) (c2)edge(d2);
			\end{scope}
			\begin{scope}[xshift=2cm]
			  \vertex (a3) at (6,2) {};
				\vertex (b3) at (6,0) {};
				\vertex (c3) at (8,0) {};
				\vertex (d3) at (8,1) {};
				\path[edge] (a3)edge(b3) (b3)edge(c3) (c3)edge(a3) (b3)edge(d3) (c3)edge(d3);
			\end{scope}
			\begin{scope}[xshift=3cm]
			  \vertex (a4) at (9,-2) {};
				\vertex (b4) at (9,0) {};
				\vertex (c4) at (11,0) {};
				\vertex (d4) at (11,-1) {}; 
				\path[edge] (a4)edge(b4) (b4)edge(c4) (c4)edge(a4) (b4)edge(d4) (c4)edge(d4);
			\end{scope}
		\end{tikzpicture}
	\end{center}
	\caption{Realizations of a graph up to rotations and translations.}
	\label{figure:realizations}
\end{figure}

Suppose that we are given a graph $G$ with edges $E$.
We consider the set of all possible realizations of the graph in the plane such 
that the lengths of the edges 
coincide with some edge labeling $\lambda\colon E\rightarrow \R_{\geq0}$. 
Edges and vertices are allowed to overlap in such a realization.
For example, suppose that $G$ has four vertices and is a complete graph minus one edge.
Figure~\ref{figure:realizations} shows all possible realizations of $G$ up to 
rotations and translations, for a given edge labeling.
We say that a property holds for a general edge labeling if it holds for all 
edge labelings belonging to a dense open subset of the vector space of all edge 
labelings. Here we address the following problem:
\begin{quote}
 \itshape{For a given graph determine the number of realizations up to
 rotations and translations for a general edge labeling.}
\end{quote}

The realizations of a graph can be considered as structures in the plane, which are 
constituted by rods connected by rotational joints.
If a graph with an edge labeling admits infinitely (finitely) many realizations 
up to rotations and translations,
then the corresponding planar structure is flexible (rigid).

The study of rigid structures, also called frameworks, was originally
motivated by mechanics and architecture, and it goes back at least to the 19th
century.
Nowadays, there is still a considerable interest in rigidity
theory~\cite{CombinatRigidity} due to various applications in natural science
and engineering.

A graph is called generically rigid (or isostatic) if a general 
edge labeling yields a rigid realization. No edge in a generically rigid graph 
can be removed without losing rigidity, that is why such graphs are also called 
minimally rigid in the literature.
The complete graph on four vertices $K_4$ is for instance not considered to be generically rigid,
since for a general choice of edge lengths it will not have a realization.
Gerard Laman~\cite{Laman1970} characterized 
the property of generic rigidity in terms of the number of edges and vertices of the graph and its 
subgraphs, hence such objects are also known as \emph{Laman graphs}.
\begin{theorem}
	\label{theorem:laman_graph}
	Let $G=(V,E)$ be a graph. 
	Then $G$ is generically rigid if and only if $\vert E\vert = 2\vert 
	V\vert-3$, and for every subgraph $G'=(V',E')$ it holds $\vert E' \vert \leq
	2\vert V' \vert-3$.
\end{theorem}

All realizations of a Laman graph can
be obtained as the solution set of a system of algebraic equations,
where the edge labels are interpreted as parameters. 
In general it is difficult to produce results on the number of real solutions of such systems.
In such situations, one often switches to a complex setting;
this also enables us to apply results from algebraic geometry.
Hence, from now on, we are interested in the number of 
complex solutions, up to an equivalence relation on $\C^2$ generalizing 
direct isometries of $\R^2$; this number is the same for any general choice of parameters, 
so we call it the \emph{Laman number} of the graph.
For some graphs up to $8$ 
vertices, this number had been computed using random values for the 
parameters~\cite{JacksonOwen2012} --- this means that it is very likely, but not 
absolutely certain, that these computations give the true numbers. Upper and 
lower bounds of Laman numbers are considered in \cite{Emiris2013,Steffens2010,Borcea2004}.
Note that for many Laman graphs there exists a labeling such that the number of real realizations is precisely the Laman number.
However, there are Laman graphs for which the Laman number gives only an upper bound on the number of real realizations.

Our main result is a combinatorial algorithm that 
computes the number of complex realizations of any given Laman graph;
it is much more efficient than just solving the 
corresponding nonlinear system of equations.
The algorithm and its correctness proof are presented in detail in \cite{SymbolicGroup}.
The purpose of the current paper is to provide a concise summary of this result focusing on the main ideas and the algorithmic point of view rather than on proofs.

\section{Main result}\label{laman_graphs}
By a \emph{graph} we mean a finite, connected, undirected 
graph without self-loops or multiple edges. We write $G = (V,E)$ to denote 
a graph $G$ with vertices $V$ and edges $E$. An edge~$e$ 
between two vertices $u$ and $v$ is denoted by $\{u,v\}$.

Using nonnegative real labels, the number of embeddings is not well-defined, since it
may depend on the actual labeling and not only on the graph.
In order to define a number that depends only on the graph, we use a
complex setting.
\begin{definition}
	\label{definition:realization}
	Let $G = (V,E)$ be a graph.
	\begin{itemize}
	  \item A \emph{labeling} of $G$ is a function $\lambda \colon E \longrightarrow \C$.
					The pair $(G, \lambda)$ is called a \emph{labeled graph}.
	  \item A \emph{realization} of $G$ is a function
					$\rho\colon V \longrightarrow\C^2$. 
					Let $\lambda$ be a labeling of $G$: we say that a realization $\rho$ is 
					\emph{compatible with} $\lambda$ if for each $e \in E$ the distance
					between its endpoints agrees with its label:
					\begin{equation}\label{equation:compatibility}
						\lambda(e) \, = \, \bigl\langle \rho(u)-\rho(v), \rho(u)-\rho(v) \bigr\rangle,\quad  e = \{u,v\},
					\end{equation}
					where $\left\langle x, y \right\rangle = x_1y_1 + x_2y_2$.
	\end{itemize}	
	A labeled graph $(G, \lambda)$ is \emph{realizable} if and only if there 
	exists a realization $\rho$ that is compatible with $\lambda$.
	
	We say that two realizations of a graph $G$ are 
	\emph{equivalent} if and only if there exists a direct
	isometry $\sigma$ of $\C^2$ between them,
	where $\sigma$ is a map of the form
	\begin{align*}
	  \vecxy
	  &\longmapsto
	  A\cdot
	  \vecxy
	  +b,
	\end{align*}
	where $A\in\C^{2\times 2}$ is an orthogonal matrix with determinant $1$ and $b\in\C^2$.
\end{definition}

\begin{definition}
	\label{definition:rigid_graph}
	A labeled graph $(G, \lambda)$ is called \emph{rigid} if it
	is realizable and	there are only finitely many realizations compatible with $\lambda$, up to 	equivalence.
\end{definition}

Our main interest is to count the number of realizations of generically rigid
graphs, namely graphs for which almost all realizable labelings induce rigidity.

\begin{definition}
	\label{definition:generically_rigid}
	A graph $G$ is called \emph{generically realizable} if for a general 
	labeling $\lambda$ the labeled graph $(G, \lambda)$ is realizable. A graph $G$ 
	is called \emph{generically rigid} if for a general labeling $\lambda$ the 
	labeled graph $(G, \lambda)$ is rigid.
\end{definition}

The number of realizations can be found by solving the system of equations
\begin{equation*}
  \bigl( (x_u - x_v)^2 + (y_u - y_v)^2 \, = \, \lambda_{uv} \bigr)_{\{u,v\} \in E}\,.
\end{equation*}
Equivalently, we can study the map $r_G$ whose preimages correspond to the solutions of the above system
\begin{equation*}
  r_G\colon \C^{2|V|} \longrightarrow \C^{E}, \quad
  (x_v, y_v)_{v \in V} \; \longmapsto \; \bigl( (x_u - x_v)^2 + (y_u - y_v)^2 \bigr)_{\{u,v\} \in E}\,.
\end{equation*}
Still we get infinitely many solutions due to translations and rotations.
Translations can be eliminated by moving one vertex to the origin.
In order to handle rotations
we perform the following transformation
$x_v\rightarrow x_v+i y_v$,\;
$y_v\rightarrow x_v-i y_v$.
Then the above equations become
\begin{align*}
   \bigl( (x_u - x_v)(y_u - y_v) \, = \, \lambda_{uv} \bigr)_{\{u,v\} \in E}\,.
\end{align*}
In this way, solutions that differ only by a rotation are the same in a suitable projective setting.
If we transform the map $r_G$ accordingly we obtain a map whose degree is finite and gives the sought number of realizations.

In order to set up a recursive formula for the degree of the map, we want to be able to handle
the two factors $(x_u - x_v)$ and $(y_u - y_v)$ independently.
To do this we duplicate the graph, and for technical reasons
we allow a more general class of graphs.
The resulting concept is roughly speaking a pair of graphs $(G,H)$ with a bijection between their sets of edges.
We identify edges by this bijection.
\begin{definition}
	\label{definition:bigraph}
	A \emph{bigraph} is a pair of finite undirected graphs $(G, H)$ ---
	allowing several components, multiple edges and self-loops ---
	where $G = (V, \edg)$ and $H = (W, \edg)$.
	The set $\edg$ is called the set of \emph{biedges}, and
	there are two maps that assign to each $e\in\edg$ the corresponding vertices in $V$ and $W$, respectively.
	Note that $G$ and $H$ are in general different graphs but there is a bijection between their sets of edges.
\end{definition}

We define the Laman number $\Lam(B)$ of a bigraph $B$ as the degree of an associated map defined in a similar way as $r_G$.
Moreover, we show (Proposition~\ref{prop:laman_number}) that the Laman number of a graph equals the Laman number of the corresponding bigraph.

The idea for proving our recursion formula is inspired by tropical geometry: we consider the 
equation system over the field of Puiseux series;
an algebraic relation between Puiseux series implies
a piecewise linear relation between their orders.
We encode these piecewise linear relations in a combinatorial data which we call \emph{bidistance}.
A bidistance is a pair of functions from the edges to rational numbers, which satisfies certain conditions.
Using a bidistance $d$ of a bigraph $B$ we can define a new bigraph $B_d$ with the same number of edges.
The solutions of the equations for the bigraph $B$ that correspond to the bidistance $d$ are in bijection with the solutions of the equations for $B_d$.
Then the solutions for $B$ are partitioned by the bidistances so that we get the following formula for the Laman number:
\begin{equation*}
  \Lam(B) \, = \, \sum_{d} \Lam(B_d).
\end{equation*}
From this we finally show the combinatorial recursion formula (Theorem~\ref{theorem:laman_number}).
For doing so we prove that $\Lam(B_d)$ is either easy to compute or the product of two Laman numbers of bigraphs with fewer edges each. 
We need some more notation to state the theorem.

\begin{definition}
	\label{definition:pseudo_laman}
	Let $B = (G, H)$ be a bigraph with biedges $\edg$, then we say that $B$ is \emph{pseudo-Laman} if
	$
	\dim(G) + \dim(H) \, = \, |\edg| + 1
	$,
	where
	$\dim(G) \, := \, |V| - |\{ \text{connected components of } G \}|$.
\end{definition}
It can be easily seen, that if $G$ is a Laman graph, then the bigraph $(G,G)$ is pseudo-Laman.
From a given bigraph we want to construct new ones with a smaller number of edges.
We introduce two constructions, a quotient and a complement, both for usual graphs and for bigraphs.
\begin{definition}
	\label{definition:quotient_graphs}
	Let $G = (V, E) $ be a graph, and let $E' \subseteq E$. We define two 
	new graphs, denoted $G \quotient E'$ and $G \subtract E'$, as follows.
	\begin{itemize}
	\item	Let $G'$ be the subgraph of $G$ determined by $E'$.
				Then we define $G	\quotient E'$ to be the graph obtained as follows:
				its vertices are the equivalence classes of the 
				vertices of $G$ modulo the relation dictating that two 
				vertices $u$ and $v$ are equivalent if there exists a path in $G'$ 
				connecting them; its edges are determined by edges in $E \setminus E'$.
	\item Let $\hat{V}$ be the set of vertices of $G$ that are endpoints of some
				edge not in $E'$. Set $\hat{E} = E \setminus E'$. Define $G \subtract E' = 
				(\hat{V}, \hat{E})$.
	\end{itemize}
\end{definition}

\begin{definition}
	Let $B = \bigl( G, H \bigr)$ be a bigraph, where $G = (V, \edg)$ and $H = (W, \edg)$.
	Given $\mcal{M} \subseteq \edg$, we define two bigraphs 
	$\leftquot{\mcal{M}}{B}  = \bigl( G \quotient \mcal{M}, \, H \subtract \mcal{M} \bigr)$ and 
	$\rightquot{\mcal{M}}{B} = \bigl( G \subtract \mcal{M}, \, H \quotient \mcal{M} \bigr)$,
	with the same set of biedges $\edg' = \edg \setminus \mcal{M}$.
\end{definition}

Finally, the main result can be stated.
\begin{theorem}
	\label{theorem:laman_number}
	Let $B=(G,H)$ be a pseudo-Laman bigraph with biedges $\edg$. Let $\bar{e}$ be a fixed biedge, then
	\begin{itemize}
		\item If $G$ or $H$ has a self-loop, then $\Lam(B) = 0$.
		\item If both $G$ and $H$ consist of a single edge joining two different vertices, then $\Lam(B) = 1$.
		\item Otherwise 
					\begin{equation}\label{equation:laman_number}
						\Lam(B) = 
						\Lam \bigl( {}^{\{ \bar{e} \} } \! B \bigr) +  
						\Lam \bigl( B^{\{ \bar{e} \}} \bigr) + 
						\sum_{(\mcal{M}, \mcal{N})} 
						\Lam \bigl( \leftquot{\mcal{M}}{B} \bigr) \cdot 
						\Lam \bigl( \rightquot{\mcal{N}}{B} \bigr),
				\end{equation}
				where each pair $(\mcal{M}, \mcal{N}) \subseteq \edg^2$ satisfies:
				\begin{itemize}
					\item $\mcal{M} \cup \mcal{N} = \edg$;
					\item	$\mcal{M} \cap \mcal{N} = \{ \bar{e} \}$;
					\item	$|\mcal{M}| \geq 2$ and $|\mcal{N}| \geq 2$;
					\item	both $\leftquot{\mcal{M}}{B}$ and $\rightquot{\mcal{N}}{B}$ are pseudo-Laman.
				\end{itemize}
	\end{itemize}
\end{theorem}
Furthermore, we show the following proposition in \cite{SymbolicGroup}, which completes the recursive algorithm for computing Laman numbers.
\begin{proposition}\label{prop:laman_number}
  The number of realizations of a Laman graph $G$ is equal to the Laman number of the bigraph $(G,G)$.
\end{proposition}

\section{Conclusion}
Using the recursion formula from Theorem~\ref{theorem:laman_number} we were able to compute Laman numbers for all Laman graphs up to 12 vertices.
Furthermore, we computed Laman numbers for single graphs up to 18 vertices, which was out of reach with the previous methods.
For further details and proofs we refer to the aforementioned paper \cite{SymbolicGroup}.
Additional information including implementations in Mathematica and in \verb!C++! can be found at \url{www.koutschan.de/data/laman/}.

\end{document}